# Boundary knot method: A meshless, exponential convergence, integration-free, and boundary-only RBF technique

W. Chen

Department of Mechanical System Engineering, Shinshu University, Wakasato 4-17-1, Nagano City, Nagano, JAPAN (E-mail: chenw@homer.shinshu-u.ac.jp)

## Abstract

Based on the radial basis function (RBF), non-singular general solution and dual reciprocity principle (DRM), this paper presents an inheretnly meshless, exponential convergence, integration-free, boundary-only collocation techniques for numerical solution of general partial differential equation systems. The basic ideas behind this methodology are very mathematically simple and generally effective. The RBFs are used in this study to approximate the inhomogeneous terms of system equations in terms of the DRM, while non-singular general solution leads to a boundary-only RBF formulation. The present method is named as the boundary knot method (BKM) to differentiate it from the other numerical techniques. In particular, due to the use of non-singular general solutions rather than singular fundamental solutions, the BKM is different from the method of fundamental solution in that the former does no need to introduce the artificial boundary and results in the symmetric system equations under certain conditions. It is also found that the BKM can solve nonlinear partial differential equations one-step without iteration if only boundary knots are used. The efficiency and utility of this new technique are validated through some typical numerical examples. Some promising developments of the BKM are also discussed.

*Key words*: boundary knot method, dual reciprocity method, BEM, method of fundamental solution,

radial basis function

# 1. Introduction

It has long been claimed that the boundary element method (BEM) is a viable alternative to the finite element method (FEM) and finite difference method (FDM). Its advantages over these domain-type methods lie mainly in dimensional reducibility and suitability to infinite domain problems. However, the practical situations are that the FEM and FDM are still dominant techniques now used in science and engineering computations. The major bottlenecks in performing the BEM have long been its weakness in handling inhomogeneous terms such as time-dependent and nonlinear problems. The recent introduction of the dual reciprocity BEM (DRBEM) by Nardini and Brebbia [1] eliminates or at least eases to great extent these inefficiencies in the BEM. However, as was pointed out in [2], the method is still more mathematically complex and requires strenuous manual effort compared with the FEM and FDM. In particular, the handling of singular integration is not easy to non-expert users and is often computationally expensive. The use of the low order polynomial approximation in the BEM also causes the slow convergence. More importantly, just like the FEM, surface mesh (or re-mesh) requires costly computation, especially for moving boundary and nonlinear problems. The method of fundamental solution (MFS) is shown an emerging technique to alleviate these drawbacks. The merits in the MFS are non-integration (of course, no singular integration), spectral convergence, and meshless (meshfree). The method is getting increasing attraction especially due to some recent works by Golberg, Chen and Muleskov [2-5]. Its shortcoming is severe ill-conditioning of the MFS equation, especially for complex boundary geometry, due to the use of artificial boundary outside physical domain. For some relative details see Kitigawa [6,7]. In addition, the present author noted that the MFS loses the symmetric matrix structure even for self-adjoint operators subject to the single type of boundary condition. These inefficiencies of the MFS motivate us to find an alternative technique, which keeps all merits of the MFS but circumvents its drawbacks.

Recently, Golberg et al. [4] established the DRBEM on the firm mathematical theory of the radial basis function (RBF). The DRBEM can be regarded as a two-step methodology. In terms of dual reciprocity method (DRM), the RBF is applied at first to approximate the particular solution of inhomogeneous terms, and then the standard BEM is used to discretize the remaining homogeneous equation. Chen et al [3] and Golberg et al [2,4] extended this RBF approximation of particular solution to the MFS, which greatly improves its applicability. In fact, the MFS itself can also be considered a special RBF collocation approach, where the fundamental solution of the governing equation is taken as the radial basis function. From the above discussions, one can conclude that the RBF is a centerpiece of the DRBEM and MFS. Therefore, there is necessity to introduce more on the RBF approach before involving the innovative work being presented in this paper.

Besides the boundary-type RBF technique of the MFS, the domain-type RBF collocation is also now under intense study since Kansa [8] first introduced the RBF to directly solve various types of partial differential system equations in 1990. The most impressive charming in the RBF-type technique is its truly meshless feature. The construction of mesh in high dimension is not a trivial work. The meshless feature is very desirable for numerical computation. Unlike the recently-developed meshless FEM with moving least square, the RBF approach is inherently simple meshless technique without any difficulty applying boundary conditions. On the other hand, Kansa and Hon [9] discovered that the convergence speed of the RBF approximation is dimensionally increased. Namely, the order of convergence rate grows linearly as the dimension increases. Therefore, the RBF avoids the dimension curse to other numerical techniques and requires much relatively fewer nodes to achieve accurate solutions in higher dimension cases. The meshless merit and dimension increasing of convergent order make the RBF approach very suitable for high-dimensional free boundary problems such as American option pricing [10]. It is obvious that the RBF is an essential factor to construct a viable numerical technique.

The purpose of this study is to introduce a novel boundary knot method (BKM) to eliminate the above-mentioned disadvantages in the MFS. The crucial point in the BKM is to use the non-singular general solution instead of the singular fundamental solution in the MFS. In this way, the inherent inefficiency of the MFS due to the use of the fictitious boundary is alleviated handily. This leads to tremendous improvement in computational stability and produces the symmetric matrix structure under certain conditions. The preliminary study of this paper shows that the BKM is a very competitive alternative to both the BEM including DRBEM and MFS. The new method can be defined as a combination of the DRM, RBF and no-singular general solution. The name of the BKM is due to its inherent meshless property, namely, the BKM does not need any discretization for any dimension problem and only uses knot point. An amazing finding of this study is the truly boundary-only technique such as the BKM may have capability to develop a non-iteration solution strategy to easily and efficiently compute the problems whose governing equations are nonlinear but boundary conditions are linear. Also time-dependent non-singular general solutions are derived to avoid the calculation of the temproal derivative for time-dependent problems.

This paper is organized as follows. Section 2 involves the procedure of the RBF approximation to particular solution. In particular, different from directive approach given in Muleskov et al [5], we use a very simple reverse method to get the analytical approximate particular solution for complicated operators and RBFs. In section 3, we introduce non-singular general solution and derive the resulting algebraic equations of the BKM. Numerical experiments and discussions are provided in section 4. Some promising potentials of the BKM are discussed in section 5, which include non-iteration nonlinear computation, non-singular general solution of some often-encountered 2D and 3D operators, time-dependent problem, localization, and conservation of symmetric matrix structure of the BKM. Finally, some concluding remarks are given based on the results reported here.

## 2. RBF approximation to particular solution

Like the DRBEM and MFS, the BKM can be viewed as a two-step numerical scheme, namely, RBF approximation to particular solution and the evaluation of homogeneous solution. The latter is the emphasis of this paper, which will be detailed in section 3. The former has been well developed by recent works of Chen and Golberg [2-4]. However, their approach [5] to determine analytical approximate solutions is rather complicated and not suitable for general case. By analogy with ideas in [11,12], we propose a viable alternative scheme to eliminate this problem very easily.

In the following, we first outline the basic methodology to employ the RBF to approximate particular solution. Let us consider the differential equation

$$L\{u(x,t)\} = f(x), \qquad x \in \Omega \tag{1}$$

with boundary conditions

$$u(x) = b_1(x), \qquad x \subset \Gamma_u, \tag{2}$$

$$\frac{\partial u(x)}{\partial n} = b_2(x), \qquad x \subset \Gamma_T, \tag{3}$$

where $L$ is a differential operator, $f(x)$ is a known forcing function. $n$ is the unit outward normal. The above simple forms of boundary conditions are to simplify description. More complicated boundary condition matters little in the present methodology. $x \in R^d$, $d$ is the dimension of geometry domain, which is bounded by a piece-wise smooth boundary $\Gamma = \Gamma_u + \Gamma_T$. It is assumed here that operator $L$ includes Laplace operator, namely,

$$L\{u\} = \nabla^2 u + L_1\{u\}. \tag{4}$$

In section 5.2, we will point out that this is not necessary. Eq. (1) can be restated as

$$\nabla^2 u + u = f(x) + u - L_1\{u\}. \tag{5}$$

The solution of the above equation (5) can be expressed as

$$u = v + u_p, \tag{6}$$

where $v$ and $u_p$ are the general (homogeneous) and particular (inhomogeneous) solutions, respectively. The latter satisfies equation

$$\nabla^2 u_p + u_p = f(x) + u - L_1\{u\}, \tag{7}$$

but does not necessarily satisfy boundary conditions (2) and (3). $v$ is the homogeneous solution of the Helmholtz equation

$$\nabla^2 v + v = 0, \quad x \in \Omega, \tag{8}$$

$$v(x) = b_1(x) - u_p, \quad x \subset \Gamma_u, \tag{9}$$

$$\frac{\partial u(x)}{\partial n} = b_2(x) - \frac{\partial u_p(x)}{\partial n}, \quad x \subset \Gamma_T. \tag{10}$$

The first step in the BKM is to evaluate the particular solution $u_p$ by the DRM and RBF. After this, Eq. (8)-(10) can be solved by the boundary collocation RBF methodology using non-singular general solution, which will be detailed in section 3.

Unless the left side of Eq. (7) is rather simple, it is practically impossible to get analytical particular solution particularly in case of nonlinear and time-dependent problems. In addition, even if the analytical solutions for some problems are available, the form of them are usually rather complicated and are not easy to use. Therefore, we prefer to approximate these inhomogeneous terms by a general numerical approach. Among them, the quasi-Monte Carlo method [13] and dual reciprocity method and RBF [2-5] are two approaches applied successfully. In this paper, we use the latter.

The DRM analogizes the particular solution by the use of a series of approximate particular solution at all specified nodes of boundary and domain. The left side of Eq. (7) is approximated by the RBF approach, namely,

$$f(x) + u - L_1\{u\} \cong \sum_{j=1}^{N+L} \alpha_j \phi(\|x - x_j\|) + \sum_{k=1}^{d} \alpha_{N+L+k} x^{(k)} + \alpha_{N+L+d+1}, \tag{11}$$

where $\alpha_j$ are the unknown coefficients. $N$ and $L$ are respectively the numbers of knots on boundary and domain. $d$ denotes the dimension of the geometry. It is noted that the use of the nodes inside domain is not necessary as will discussed in later section 5.1. $\|\bullet\|$ represents the Euclidean norm, $\phi(\|x - x_j\|)$ is the RBF function of one-dimensional distance variable. The first-order polynomial function in Eq. (11) is set up to assure nonsingularity of interpolation matrix, since most of the globally support RBFs are only conditionally positive definite [14]. Note that $x^{(k)}$ indicates different independent variables.

By forcing Eq. (11) to exactly satisfy Eq. (7) at all nodes together with the following $d+1$ side conditions

$$\sum_{j=1}^{N+L} \alpha_{N+L+j} x_j + \alpha_{N+L+d+1} = 0, \tag{12}$$

we can uniquely determine

$$\alpha = A_\phi^{-1} \begin{Bmatrix} f(x) + u - L_1\{u\} \\ 0 \end{Bmatrix}, \tag{13}$$

where $A_\phi$ means nonsingular interpolation matrix and holds symmetric structure. Similar to interpolation formula (11), applying the RBF interpolation to approximation of function $u$ yields

$$u \cong \sum_{j=1}^{N+L} \lambda_j \phi(\|x - x_j\|) + \sum_{k=1}^{d} \lambda_{N+L+k} x^{(k)} + \lambda_{N+L+d+1}, \tag{14}$$

namely,

$$\begin{Bmatrix} u \\ 0 \end{Bmatrix} = A_\phi \lambda. \tag{15}$$

For operator $L_1$, we have

$$L_1 \begin{Bmatrix} u \\ 0 \end{Bmatrix} = L_1 \{A_\phi\} \lambda. \tag{16}$$

So

$$L_1 \begin{Bmatrix} u \\ 0 \end{Bmatrix} = L_1 \{A_\phi\} A_\phi^{-1} \begin{Bmatrix} u \\ 0 \end{Bmatrix}. \tag{17}$$

Substituting Eq. (17) into Eq. (13), we get

$$\alpha = A_\phi^{-1} \left[ \begin{Bmatrix} f(x) \\ 0 \end{Bmatrix} + \left(I - L_1\{A_\phi\} A_\phi^{-1}\right) \begin{Bmatrix} u \\ 0 \end{Bmatrix} \right]. \tag{18}$$

In the normal procedure, the RBF $\phi$ in Eqs. (11) and (14) is first selected, and then corresponding approximate particular solutions $\varphi$ are determined by analytically integrating equation

$$\nabla^2 \varphi + \varphi = \begin{Bmatrix} \phi \\ x^{(k)} \\ 1 \end{Bmatrix}, \qquad j = 1, \mathrm{K}, N+L,, \quad k = 1, \mathrm{K}, d. \tag{19}$$

This methodology performs easily only for simple Laplace operator and linear RBF

$$\phi(x, x_j) = \phi(r_j) = r_j, \tag{20}$$

where $r_j = \|x - x_j\|$. The other more popular RBFs include multiquadratics (MQ)

$$\phi(x, x_j) = \phi(r_j) = \sqrt{r_j^2 + c_j^2}, \tag{21}$$

where $c_j$ is the shape parameter, and thin plate spline (TPS)

$$\phi(x, x_j) = \phi(r_j) = r_j^{2m} \log(r_j), \tag{22}$$

where $m$ is order of TPS. Muleskov et al. [5] recently got the TPS analytical approximate particular solutions of Helmholtz operator by using a sophisticated mathematical approach, and the forms of their solutions are rather complicated. In addition, the analytical approximate particular solutions for general cases such as the MQ and other differential operators are not yet available now due to great difficulty involved.

In this study, we use an alternative reverse approach to eliminate this difficulty in a much easier fashion. Inverse to the normal scheme, the proposed approach is to determine the approximate particular solution $\varphi$ at first, and then we can very easily evaluate the corresponding RBF $\phi$ through substituting the known particular solution into Eq. (19). The principle used here to choose approximate solution for the Helmholtz operator is to increase two-order power of the often-used RBFs and to avoid the singularity at $r=0$. For the MQ, second-order TPS ($m=2$) and linear RBFs, the chosen approximate particular solution are respectively

$$\varphi(r_j) = (r_j^2 + c_j^2)^{3/2}, \tag{23a}$$

$$\varphi(r_j) = r_j^6 \log(r_j), \tag{23b}$$

$$\varphi(r_j) = r_j^3. \tag{23c}$$

The corresponding RBFs are

$$\phi(r_j) = 6(r_j^2 + c_j^2) + \frac{3r^2}{\sqrt{r_j^2 + c_j^2}} + (r_j^2 + c_j^2)^{3/2}, \tag{24a}$$

$$\phi(r_j) = 36 r_j^4 \log(r_j) + 12 r_j^4 + r_j^6 \log(r_j), \tag{24b}$$

$$\phi(r_j) = 9 r_j + r_j^3. \tag{24c}$$

The present strategy can also be easily extended to the 2D, 3D biharmonic or convection diffusion operators other than the present 2D Helmholtz operator. The numerical experiments given in section 4 demonstrate that this scheme works well and robustly. For example, we observe of the super convergence of the MQ clearly.

Finally, we can get particular solutions at any point by summation of approximate particular solutions at all nodes, namely,

$$u_p = \sum_{j=1}^{M} \alpha_j \varphi(\|x - x_j\|), \qquad (25)$$

where $M = N + L + d + 1$. Substituting Eq. (18) into Eq. (25) yields

$$u_p = \Phi A_\phi^{-1} \left[ \begin{Bmatrix} f(x) \\ 0 \end{Bmatrix} + \left(I - L_1\{A_\phi\}A_\phi^{-1}\right)\begin{Bmatrix} u \\ 0 \end{Bmatrix} \right]. \qquad (26)$$

Note that $\Phi$ here is known matrix comprised of $\varphi(\|x - x_j\|)$. Eq. (26) will be used in the boundary collocation equations introduced in next section 3.

## 3. Non-singular general solution and boundary knot method

It is well known that the Laplace operator of the simplest among all the known operators has not nonsingular general solution. Golberg and Chen [2] pointed out that the placement of source points outside domain in the MFS is to avoid the singularities. However, we found through numerical experiments that even if the source and response points on physical boundary are chosen differently without the singularities, the solutions were still degraded evidently compared with those by the standard MFS, where the source points are placed only on the fictitious auxiliary boundary. Moreover, the more distant the source points are located from physical domain, the more accurate solutions are obtained for the MFS [2]. However, unfortunately, the resulting equations becomes extremely ill-conditioned which in some cases deteriorates the solution [2,6,7].

To illustrate the basic idea of the boundary collocation using non-singular general solution, we take the 2D Helmholtz operator as an illustrative example. There do exist non-singular general solutions for the other differential operators. The reason for this choice is that the Helmholtz operator is the simplest among various often-encountered operators having non-singular general solution.

The 2D homogeneous Helmholtz equation (8) has two general solutions, namely,

$$v(r) = c_1 J_0(r) + c_2 Y_0(r), \tag{27}$$

where $J_0(r)$ and $Y_0(r)$ are the zero-order Bessel functions of the first and second kinds, respectively. In the standard BEM and MFS, the Hankel function

$$H(r) = J_0(r) + i Y_0(r) \tag{28}$$

is applied as the fundamental solution. It is noted that $Y_0(r)$ encounters logarithm singularity, which causes the major technique difficulty in applying the BEM and MFS. Many special techniques have developed to solve or circumvent this troublesome singular problem.

The crucial process in the present BKM scheme is to discard this singular $Y_0(r)$ general solution and to only use $J_0(r)$ as the radial basis function to collocate the boundary condition equations (9) and (10). It is noted that $J_0(r)$ exactly satisfies the Helmholtz equation and we can therefore get a boundary-only RBF collocation scheme. Unlike the MFS, all collocation nodes are placed only on physical boundary and can be used as either source or response points. The previously-mentioned drawbacks in the MFS disappear naturally in the present BKM.

Let $\{x_k\}_{k=1}^{N}$ denote a set of nodes on the physical boundary, the solution $v(x)$ of Eq. (8) is approximated in a standard RBF collocation procedure

$$v(x) = \sum_{k=1}^{N} \beta_k J_0(r_k), \tag{29}$$

where $r_k = \|x - x_k\|$. $k$ is index of source points. $N$ is the total number of boundary knots. $\beta_k$ are the desired coefficients.

In terms of the collocation of Eqs. (9) and (10), we have

$$\sum_{k=1}^{N} \beta_k J_0(r_{ik}) = b_1(x_i) - u_p(x_i), \tag{30}$$

$$\sum_{k=1}^{N} \beta_k \frac{\partial J_0(r_{jk})}{\partial n} = b_2(x_j) - \frac{\partial u_p(x_j)}{\partial n}, \tag{31}$$

where $i$ and $j$ indicate respectively corresponding boundary response knots. If internal node are used, we need constitute another set of supplement equations. It is evident that the following equations at interior knots

$$\sum_{k=1}^{N} \beta_k J_0(r_{lk}) = u_l - u_p(x_l), \quad l = 1, \mathrm{K}, L \tag{32}$$

can be established, where $l$ indicates the internal response knots and $L$ is the total number of interior points used. Substituting Eq. (26) into Eqs. (30), (31) and (32), we get $N+L$ simultaneous algebraic equations. The unknown variables of these equations are boundary collocation coefficients $\beta_k$ and values of domain-inside dependent variable $u_l$. After solving these equations, we need compute coefficient $\alpha$ by Eq. (18) so that we can get the approximate solution at any points by

$$u(x) = \sum_{k=1}^{N} \beta_k J_0(\|x - x_k\|) + u_p, \tag{33}$$

where $u_p$ is evaluated by Eq. (26).

It is stressed that the use of interior points are not necessary in the BKM. The term boundary-only is used here in the sense as in the DRBEM and MFS that only boundary knots are required, although the use of internal knots can improve solution accuracy. If only boundary knots are employed, Eq. (32) should be omitted and $u_p$ and $\alpha$ are determined directly by Eq. (18) and (26) before the solution of resulting algebraic Eqs. (30) and (31). From the above analysis, one can find that the present BKM is rather mathematically simple and easy-to-program. The extension of this methodology to the other partial differential operators is very straightforward.

Before proceeding the numerical experiments, we consider to choose the radial basis function. In general, RBFs are globally supported basis functions and lead to a dense matrix, which becomes highly ill-conditioned if very smooth radial basis functions are used with large number of interpolation nodes [14]. This causes severe stability problems and computationally inefficiency for large size problem. Some approaches have proposed to remedy this problem such as domain decomposition and compactly supported RBFs (CS-RBFs). The latter is very recently developed by Wendland [15], Wu [16] and Schaback [17]. The CS-RBFs can result in sparse banded matrices and is therefore preferred in many cases. Golberg et al. [2] and Wong et al. [18] respectively applied the CS-RBFs to the MFS and RBF collocation method successfully. An excellent comprehensive research on this topic can be found in Kansa and Hon [9]. Chen et al [19] recently also used the CS-RBFs to the DRBEM. However, in this study, we will do not use the CS-RBFs to focus on the illustration of the basic idea of the BKM with simpler globally-supported RBFs. Possible benefits of the CS-RBFs to the BKM will be discussed in later section 5.4.

The MQ, TPS and linear RBF are very simple and most widely used now, which were introduced respectively by Hardy [20], Duchon [21], and Nardini and Brebbia [1]. Among them, it is well known that the MQ ranks the best in accuracy [22] and is only the RBF with desirable merit of spectral convergence [4]. However, its accuracy is greatly influenced by the choice of the shape parameter [23,24]. So far the optimal choice of shape parameter is still an open research topic. In spite of this problem, the MQ is still widely used in the RBF solution of various differential systems. For numerical example of Laplace equation shown in next section 4.2, the linear and TPS RBFs have slower convergence rate than the MQ. For example, the average relative error is 0.91% for the Linear RBF with 11 knots, -0.39% for the TPS with 11 knots and –0.023% for the MQ (shape parameter 2) with 9 knots. We even got the very accurate solutions, 0.50% relative average error, by the MQ with only 3 knots. Also the use of the other type of RBFs such as the TPS is satisfactory in many cases. To simplify the presentation, this paper only uses the MQ to approximate the particular solution.

## 4. Numerical results and discussions

In this paper, all numerical examples unless otherwise specified are extracted from [25]. The geometry of test problem is all an ellipse featured with semi-major axis of length 2 and semi-minor axis of length 1 as shown in Figs. 1 and 2. These examples are chosen since their analytical and numerical solutions are obtainable so that the prowess of the new BKM technique can explicitly exposed through comparison. More complex problems can be handled without any extra difficulty. The Bessel function of the first kind of zero order is evaluated via a short program given in [26]. The 2D Cartesian co-ordinates (x,y) system is used as in [25].

### 4.1. Helmholtz equation

The governing equation of the 2D Helmholtz type is given by

$$\nabla^2 u + \lambda^2 u = 0. \tag{34}$$

In this study, we let =1 and pose inhomogeneous boundary condition

$$u = \sin x. \tag{35}$$

It is obvious that Eq. (35) is also a particular solution of Eq. (34). The BKM collocation of this test example is shown in Fig. 1. The DRBEM discretization with 17 interior points is illustrated in Fig. 2. Numerical results by the present BKM is displayed in Table 1 together with those by the DRBEM for comparison.

The numbers in brackets of Table 1 mean the total nodes used. It is found that the present BKM converges very quickly. The BKM solutions using 7 nodes are also satisfactory. This demonstrates that the BKM enjoys the super-convergent property as in the other types of collocation methods [27]. In contrast, the standard BEM has only low order of convergence speed due to its lower order polynomial approximation [2]. Please note that in this case, there is no particular solution to be approximated in the BKM by using the RBF and DRM since the inhomogeneous term does not appear.

**4.2. Laplace equation**

Reader may argue that it is somehow unfair to choose the homogeneous Helmholtz equation to compare the BKM and DRBEM which uses the Laplace fundamental solution in the above example. To further justify the exponential convergence of the BKM compared with the standard BEM, in the following we exam Laplace equation

$$\nabla^2 u = 0 \tag{36}$$

with boundary condition

$$u = x + y. \tag{37}$$

Eq. (37) is easily found to be a particular solution of Eq. (36). This problem is typically well suited to be solved by the standard BEM technique. Therefore, it is a good example to show up the accuracy and efficiency of the BKM vis-a-vis the BEM. In order to apply the non-singular solution of the Helmholtz operator, Eq. (36) is rewritten as

$$\nabla^2 u + u = u. \tag{38}$$

The left inhomogeneous term $u$ is approximated by the MQ in terms of the DRM as shown in the previous section 2. The numerical results are displayed in Table 2, where the BEM solutions come from [25].

The MQ shape parameter is $c=25$ for both 3 and 5 boundary knots in the BKM. It is observed that the solutions are not sensitive to the parameter $c$ in this case. It is surprising to see from Table 2 that the BKM solutions using 3 boundary knots achieve the accuracy of four significant digits and are much more accurate than the standard BEM solution using 16 boundary nodes. This striking accuracy of the BKM with very few knots again validates its spectral convergence. No interior point is employed in both the BEM and the present BKM in this case. It is interesting to note that the coefficient matrix in the BEM and BKM are both fully populated. However, unlike the BEM, the BKM has systematic coefficient matrix for all self-adjoint operators with only one type of boundary conditions as this case. In next section 5.4, we will propose some possible strategy to conserve this symmetric matrix structure in the BKM to the problems with several types of boundary conditions.

**4.3. Convection diffusion equation**

The FEM and FDM encounter some difficulty to produce accurate solution to the systems involving the first order derivative of convection term. Special care need be taken to handle this problem. It is claimed that the BEM does not suffer these accuracy problems occurring in the FEM and FDM in such cases. In particular, the DRBEM was said to be very suitable for this type problem [25]. As a starting point applying the BKM to general system equation, we consider the convection diffusion equation

$$\nabla^2 u = -\partial u/\partial x, \tag{39}$$

which is given in [25] to test the DRBEM. The boundary condition is stated as

$$u = e^{-x}, \tag{40}$$

which also constitutes a particular solution of this problem. By adding u on dual sides of Eq. (39),

we have

$$\nabla^2 u + u = u - \partial u/\partial x. \tag{41}$$

The above Eq. (41) is suitable to be solved by the BKM. The results by both the DRBEM and BKM are listed in Table 3.

The MQ shape parameter is chosen 4 in this calculation. The BKM uses 7 boundary knots and 8 or 11 internal knots. In contrast, the DRBEM [25] applied 16 boundary and 17 inner nodes. It is noted that in this case the use of the interior points can improve the solution accuracy evidently. This is due to the fact that the governing equation is convection-diffusion domain-dominant problem. The Helmhotlz operator in the present BKM and Laplace operator in the DRBEM [25] can not capture well convection effects of the system equation only by using boundary nodes. In next section 5.1, we will give some further discussions to this type of problems. It is found from Table 3 that the BKM achieves the salient accurate solutions and outperforms the DRBEM in computational efficiency due to the super-convergent features of the MQ interpolation and global BKM collocation scheme.

Further consider equation

$$\nabla^2 u = -\partial u/\partial x - \partial u/\partial y \tag{42}$$

with boundary conditions

$$u = e^{-x} + e^{-y}. \tag{43}$$

It is obvious that the boundary condition equation (43) is also a solution of Eq. (42). The numerical results are summarized in Table 4.

In the BKM, the MQ shape parameter is taken 5.5. We employ 8 or 11 inner points and 7 boundary knots in the BKM compared with 16 boundary nodes and 17 inner points in the DRBEM [25]. It is seen that the BKM works equally in this case. It is found that the BKM with fewer points produced almost the same accurate solutions as in the DRBEM. Considering the extremely simple and easy-to-use merits in the BKM, the method is obviously superior to the DRBEM.

## 5. Further developments

In the foregoing section 4, we only dealt with the two-dimensional linear stationary problems. By its very basis, the BKM can be handily extended to the nonlinear, three-dimensional, time-dependent, and any types of partial differential systems. This section will brief some relative theoretical and numerical results we have so far finished. The purpose is to expose great promising potentials of the BKM in some important aspects.

### 5.1. Nonlinear problems

One of the major challenges in computational engineering and science is to solve the nonlinear problem. The standard solution procedure of nonlinear modeling equations is mostly originated from the Newton method and its variants [28]. For these methods, the proper determination of initial guess solution and evaluation of Jacobian matrix and its inverse are either difficult or very computationally expensive tasks. It is highly desirable to circumvent or ease to some extent these difficulties. It is noted that most of nonlinear systems are comprised of nonlinear governing equations and linear boundary conditions. This reminds us that the boundary-only technique may provide such a possible strategy to develop one-step nonlinear analysis technique. In particular, the BKM, MFS and DRBEM are feasible choice due to their truly boundary-only merit for general nonlinear problems compared with normal BEM techniques, which are often cumbersome in solving nonlinear system due to the fact that it is nearly impossible to find fundamental solution of most practically significant nonlinear problems.

The basic methodology presented here is that we should as possible not use the internal points in the BKM solution of nonlinear problems. To clearly state our idea, consider the following nonlinear equation in two dimension

$$\nabla^2 u + u_{xx} u = 2 + 2x^2 \tag{44}$$

with boundary condition

$$u = x^2, \tag{45}$$

which is also the solution of this problem. Adding $u$ to both sides of Eq. (44), we have

$$\nabla^2 u + u = u - u_{xx} u + 2 + 2x^2. \tag{46}$$

Before further modeling, we here introduce concept of the Hadmard matrix product to simplify the nonlinear formulation expression [29].

**Definition 5.1** Let matrices $A=[a_{ij}]$ and $B=[b_{ij}] \in C^{N \times M}$, the Hadamard product of matrices is defined as $A \circ B = [a_{ij} b_{ij}] \in C^{N \times M}$. where $C^{N \times M}$ denotes the set of $N \times M$ real matrices.

Following the standard BKM procedure, we have

$$\alpha = A_\phi^{-1} \left[ u - \left( \frac{\partial^2 A_\phi}{\partial x^2} A_\phi^{-1} u \right) \circ u + 2x^2 + 2 \right] \tag{47}$$

and

$$u_p = \Phi \alpha,$$

where $A_\phi$ is the RBF interpolation matrix as given in Eq. (13). If only boundary knots are employed, we will have

$$Jc = x^2 - \Phi \left\{ A_\phi^{-1} \left[ x^2 - \frac{2}{x} \circ \left( \frac{\partial^2 A_\phi}{\partial x^2} A_\phi^{-1} x^2 \right) + 2x^2 + 2 \right] \right\} = b(x), \tag{48}$$

where $J$ and $c$ are boundary interpolation matrix of the Bessel function and the corresponding unknown coefficient vector. Eq. (48) is linear algebraic equations of boundary modeling, which can be solved readily. Table 5 displays some computed results. Five boundary knots are used with the MQ shape parameter 2. The average relative error is 0.91%. It is seen from Table 5 that the accuracies of solutions at all knots are satisfactory.

However, the present scheme does not work equally well for any cases. For example, let us consider the Berger equation of steady-state situation

$$\nabla^2 u + u_x u = 0. \tag{49}$$

Here we still use data given in [25], where the boundary condition is given by

$$u = 2/x. \tag{50}$$

Eq. (50) is also a particular solution of Berger equation (49). The BKM solutions with only boundary knots are listed in Table 6.

The shape parameter is taken 1 in the MQ. The average relative absolute error is 3.97%. It is worth stressing that the Berger equation has the structure of convection diffusion equation and $u$ in Eq. (49) is in fact a velocity. Therefore, it is the convection term which causes the inaccuracy of the solution if not using inner points. We also found in the previous linear convection diffusion cases 4.3 that the solution accuracy is degraded to 10.7% if without interior points and is improved to 0.56% if with few interior points. The similar problem also happened in DRBEM solution of the domain-dominant convection diffusion problem using the Laplace fundamental solution. As was pointed out in [25], the internal nodes used in DRM analysis is not in general a condition to obtain accurate solution. In most cases, the use of some interior nodes is useful to improve the solution accuracy. The inaccuracy in the present case is mainly due to the inefficiency in the RBF approximation of particular solution with convection characteristics if without using interior nodes. If the proper RBFs are applied to capture the system convection effects accurately in approximation of particular solution, the internal points will not be necessary to improve accuracy. By analogy with multiple reciprocity BEM [30], it may be feasible to apply the zero-order Bessel function of the first kind repeatedly as the RBF to approximation of particular solution, which may be superior to the MQ in this case. We will investigate this scheme in next stage of work.

It is expected that if inner points are used in the present formulation, the solution accuracy will be improved greatly as in the DRBEM [25]. However, the use of internal node will lead to nonlinear algebraic systems of equations although it may be much smaller in magnitude than those by other methods. In the DRBEM, it is reported [25] that the use of the fundamental solution of convection diffusion equation can significantly improve the solution accuracy of linear transient convection diffusion problem without using inner point. This suggests us that a non-singular general solution of convection-diffusion equation may be much more suitable for this type of problem. For two-dimensional steady-state convection-diffusion equation

$$D\nabla^2\phi + v_x \frac{\partial \phi}{\partial x} + v_y \frac{\partial \phi}{\partial y} + k\phi = 0, \tag{51}$$

where $v_x$ and $v_y$ are the components of the velocity vector $v$, $D$ is the diffusivity coefficient and $k$ represents the reaction coefficient. We have non-singular general solution

$$\phi^* = \frac{1}{2\pi D} e^{-\frac{v \bullet r}{2D}} J_0(\mu r), \tag{52}$$

where

$$\mu = \left[\left(\frac{|v|}{2D}\right)^2 + \frac{k}{D}\right]^{\frac{1}{2}} \tag{53}$$

and marked $r$ denotes the distance vector between the source and response knots. $J_0$ is the zero-order Bessel function of the first kind. The above non-singular general solution is in fact the non-singular part of the fundamental solution of convection diffusion equation. According to the experience in the DRBEM, alternative non-singular general solution (52), which reflects the convection effects, may make the use of the interior points unnecessary in the BKM for the previous convection diffusion case 4.3 [25]. It is expected that the same thing will happen to the BKM solution of nonlinear Burger equation. It is worth stressing that since in nonlinear cases velocity values in non-singular general solution (52) may vary at different boundary knots, the RBFs used on distinct boundary response knots may differ in velocity coefficients. Namely, source knot-dependent RBFs need be used. In this way, the resulting algebraic matrix will be linear if the only boundary nodes are applied.

The convection diffusion is not an easy problem to handle. The oscillation and damping manifest in the FEM and FDM solution of this type of problems. Even if the interior knots are necessarily used to assure accurate solution in such case, the number of the knots in the BKM will be much fewer than the domain-type methods. Consequently, the BKM will results in a much less magnitude of nonlinear algebraic equations. Moreover, the results by the BKM without interior nodes can provide a reliable initial guess solution for iteration computation. The research on these interesting issues are now under way and are left to address in detail in the subsequent paper. On the other hand, it is conceivable that we can use the same strategy to handle nonlinear boundary condition problem. Namely, inverse to the present BKM scheme, we may choose the RBFs, which satisfy the nonlinear boundary condition equation. Then we have a domain-type linear RBF collocation modeling.

Recently Chen et al. [29] demonstrated that the special matrix product approach is very promising and powerful techniques to nonlinear stability analysis, uncoupling computation, evaluation of Jacobian matrix, and construction of simple, efficient and robust iteration formulas. The Hadamard product and SJT product are two types of special matrix product very useful to all point-approximation techniques, which include the finite difference, finite volume, collocation, and RBF-related techniques (DRBEM, MFS and the present BKM). In this section, we only used the Hadamard product to simplify the representation of the RBF nonlinear formulation. A significant research topic will be explore the potential in a combined use of the Hadamard and SJT product technique and the BKM to solve the nonlinear problems.

**5.2. Non-singular general solutions for some differential operators and others**

In the following, we will give non-singular general solutions of some often-used 2D and 3D operators. In the following equations $r$ means distance value and $c$ is undetermined coefficients. For 2D equation

$$\nabla^2 u - \lambda^2 u = 0, \tag{54}$$

the non-singular particular solution is

$$u^* = cI_0(\lambda r), \tag{55}$$

where $I_0$ is the zero-order modified Bessel function of the first kind. For the 3D Helmholtz-like operator,

$$\nabla^2 u \pm \lambda^2 u = 0, \tag{56}$$

we have the non-singular general solution

$$u^* = c\frac{\sin(\lambda r)}{r}, \tag{57}$$

and

$$u^* = c\frac{\sinh(\lambda r)}{r}, \tag{58}$$

where sinh denotes the hyperbolic function. For the 2D biharmonic operator

$$\nabla^4 w - \lambda^2 w = 0, \tag{59}$$

we have the non-singular general solution

$$w^* = c_1 J_0(\lambda r) + c_2 I_0(\lambda r), \tag{60}$$

where $J_0$ and $I_0$ are respectively the Bessel and modified Bessel functions of the first kind of the zero-order. The non-singular general solution of the 3D biharmonic operator is given by

$$w^* = c_0 \frac{\sin(\lambda r)}{r} + c_1 \frac{\sinh(\lambda r)}{r}, \tag{61}$$

Katsikadelis and Nerantzaki [31] proposed a simple strategy to apply the DRBEM to any partial differential equations which do not include standard Laplace or biharmonic operators. The strategy can be readily extended to the BKM calculation of any linear or nonlinear systems yet keeping its boundary-only merit.

### 5.3. Time-dependent problems

For time-dependent problem, a numerical time integrator is often used to approximate time derivative after spatial discretization is done. The use of marching scheme of time integrator involves some difficult issues relating to stability and accuracy [32]. Another methodology is the mode analysis. However, the method is not very applicable for many cases such as shock. We here hope to introduce an alternative technique free from computing temporal derivative within the BKM frame. To explain clearly, let us consider 2D heat and diffusion governing equation

$$\Delta u = \frac{1}{k} \frac{\partial u}{\partial t}. \tag{62}$$

By separation, we have

$$u = \phi(r)\varphi(t). \tag{63}$$

Substituting the above equation into governing equation (62) yields

$$\varphi \Delta \phi = \frac{1}{k} \phi \frac{\partial \varphi}{\partial t}. \qquad (64)$$

Let

$$\phi = J_0(r), \qquad (65)$$

where $J_0$ is the Bessel function of the first kind of zero order, we have

$$\varphi' + k\varphi = 0. \qquad (66)$$

Its general solution is

$$\varphi(t) = Ae^{-kt}. \qquad (67)$$

Therefore, we have general solution of transient heat and diffusion equation

$$u(r,t) = Ae^{-kt}J_0(r), \qquad (68)$$

where $A$ is unknown constants determined by boundary and initial conditions. Furthermore, consider wave equation

$$\Delta u = \frac{1}{c^2}\frac{\partial^2 u}{\partial t^2}. \qquad (69)$$

Similarly, we get the general solution

$$u(r,t) = [A\cos(ct) + B\sin(ct)]J_0(r), \qquad (70)$$

which satisfies the governing equation (69). By using the above non-singular time-dependent general solutions (68) and (70) as the RBFs, we can obtain the BKM boundary-only formulation for time-dependent diffusion problems without time derivative. The strategy is straightforward and in fact equivalent to the BEM using time-dependent fundamental solution. The difficulty in applying this scheme may lie in how to satisfy the initial conditions inside domain as in the BEM.

Recently Zhu et al. [33] used the Laplace transforms to eliminate temporal derivative in the DRBEM solution of time-dependent problems. The accuracy and efficiency of their method are satisfactory. More recently Golberg and Chen [2] combined this approach with the MFS to solve the transient Helmholtz problems successfully. The key to successful use of the Laplace transforms is to implement the Stehfest s algorithm [34] for the accurate numerical inversion of Laplace transform. It is very obvious that such Laplace transforms methodology can be extended to the BKM solution of transient problems.

### 5.4. Symmetricity, localization and system conditioning

Compared with finite elements, some evident drawbacks in the DRBEM are full matrix, loss of the symmetric matrix structure for self-adjoint operators, and singular numerical integration. Unlike the DRBEM, the BKM has not singular integration and is very easy to use due to its integration free merit. However, it is noted that the global support BKM leads to the fully-populated system matrix. For large system problem, the solving equations of full matrix of large size is rather computationally costly. Like the MFS, there are two issues relating to the full coefficient matrix in the BKM, namely, globally-supported boundary collocation using non-singular general solution and RBF approximation of particular solutions. Recently, Chen et al. [19] applied compactly-supported RBF (CS-RBF) to obtain the banded system matrix in the RBF approximation of particular solution. Obviously, this methodology is applicable to the RBF and DQM approximate particular solution in the present BKM.

On the other hand, as was mentioned previously, the BKM also loses the symmetric matrix structure for self-adjoint operators due to more than one type of boundary conditions. The same thing happens to the RBF collocation method. A Hermite interpolation RBF has been proposed to preserve the symmetric structure of the resultant coefficient matrix. The strategy is also applicable to the BKM. Namely, the boundary RBF collocation in the BKM will be performed locally within the same type of boundary condition separately. The different collocation approximate equations in local domains will then be matched through keeping $C_0$ continuity at interface knots in such a way that we get a sparse and symmetric banded system matrix of boundary collocation approximation. As was pointed out in [35], the symmetric matrix can assure the good-conditioning of system equation. The further investigation in this direction will be significant.

A combines use of the CS-RBF and local Hermite RBF collocation interpolation on the boundary builds the BKM a local numerical modeling technique as a whole. The localization of the BKM can eliminate the ill-conditioning in the global BKM when applied to large system problem with a huge number of nodes. Numerical investigation will be subject of further study along this line.

## 6. Concluding remarks

The present BKM can be regarded one kind of the Trefftz method [36], where the trial function is required to satisfy governing equation. The BKM distinguishes the other Trefftz techniques in that the RBF is used through non-singular general solution. The BKM also looks like the MFS in many senses except that the non-singular general solution is used instead of the singular fundamental solution. The major shortcoming of the MFS to use fictitious boundary is no longer inherent in the present BKM and the robustness of solution is greatly improved. On the other hand, the BKM may be especially attractive in the solution of three-dimensional problem due to the RBF s dimensionally-increased order of convergent rate recently found by Kansa and Hon [9].

The proper choice of the RBF to approximate the particular solution is one of pivotal issues affecting the BKM efficiency. Although the MQ holds exponential convergence merit, the problem-dependent shape parameter degrades its attractiveness. Among all currently available RBFs, the CS-RBF and TPS are more practically useful in general. We plan to combine these two types of RBFs with the BKM to investigate some benchmark problems in future work. In particular, the CS-RBF can lead to a banded systematic interplant matrix and is therefore preferred. On the other hand, it is still possible to develop some new RBFs which at least partly comply the intrinsic characteristics of the targeted problems as done very successfully in [12,13] for the DRBEM solution of axisymmetric and infinite diffusion problems. For example, even in the domain-type collocation RBF method, it is expected that the zero-order Bessel of the first kind and non-singular general solution of convection diffusion equation may be the RBF of the choice in the collocation of various wave, convection, heat, and diffusion problems. The same care should also be taken in the DRM and RBF approximation of particular solution within the BKM.

In conclusion, the BKM presented in this paper is a novel technique combining several powerful methodologies of the DRM, RBF, and non-singular general solution. The method inherently possesses desirable numerical merits which include meshless, boundary-only, non-integration, exponential convergence, and mathematical simplicity. The use of the method is remarkably easy, especially for complicated shapes and higher dimensions. By using the local Hermite interpolation and compactly-supported RBF, the resulting BKM approximate matrix as a whole is expected to have a symmetric sparse banded structure for self-adjoint operator problem with any boundary conditions. The possibility to develop a non-iteration BKM nonlinear solver with only boundary nodes is especially practically significant, which will be emphasis of the next recent future work. In addition, more numerical experiments to test the BKM will be beneficial. This paper can be regarded a starting point of a series of work afterward.

**References**:

1. D. Nardini and C.A. Brebbia, A new approach to free vibration analysis using boundary elements.

*Applied Mathematical Modeling*, **7** 157-162 (1983).

2. M.A. Golberg and C.S. Chen, The method of fundamental solutions for potential, Helmholtz and diffusion problems. In *Boundary Integral Methods - Numerical and Mathematical Aspects,* (Edited by M.A. Golberg), pp. 103-176, Computational Mechanics Publications, (1998).

3. C.S. Chen, The method of potential for nonlinear thermal explosion, *Commun Numer. Methods Engng*, **11** 675-681 (1995).

4. M.A. Golberg, C.S. Chen, H. Bowman and H. Power, Some comments on the use of radial basis functions in the dual reciprocity method, *Comput. Mech.* **21** 141-148 (1998).

5. A.S. Muleskov, M.A. Golberg and C.S. Chen, Particular solutions of Helmholtz-type operators using higher order polyharmonic splines, *Comput. Mech.* **23** 411-419 (1999).

6. T. Kitagawa, On the numerical stability of the method of fundamental solutions applied to the Dirichlet problem. *Japan Journal of Applied Mathematics*, **35** 507-518, (1988).

7. T. Kitagawa, Asymptotic stability of the fundamental solution method. *Journal of Computational and Applied Mathematics*, **38** 263-269 (1991).

8. E.J. Kansa, Multiquadrics: A scattered data approximation scheme with applications to computational fluid-dynamics. *Comput. Math. Appl*. **19** 147-161 (1990).

9. E.J. Kansa and Y.C. Hon, Circumventing the ill-conditioning problem with multiquadric radial basis functions: applications to elliptic partial differential equations. *Comput. Math. Appls*. **39** 123-137 (2000).

10. Y.C. Hon and X.Z. Mao, A radial basis function method for solving options pricing model. *Financial Engineering*, **81**(1) 31-49 (1999).

11. L.C. Wrobel, J.C.F. Telles and C.A. Brebbia, A dual reciprocity boundary element formulation for axisymmetric diffusion problems. In *Boundary Element VIII*, Vol. **1**, Comput. Mech. Publ., Southampton, UK (1986).

12. C.F. Loeffler and W.J. Mansur, Dual reciprocity boundary element formulation for potential problems in infinite domains. In *Boundary Elements X*, Vol. **2**, Comput. Mech. Publ. and Springer-Verlag, Berlin (1988).

13. C.S. Chen, M.A. Golberg and Y.C. Hon, The method of fundamental solutions and quasi-monte-carlo method for diffusion equations. *Int. J. Numer. Meth. Engng*. **43** 1421-1435 (1998).

14. M. Zerroukat, H. Power and C.S. Chen, A numerical method for heat transfer problems using collocation and radial basis function, *Inter. J. Numer. Method Engng*. **42** 1263-1278 (1998).

15. H. Wendland, Piecewise polynomial, positive definite and compactly supported radial function of minimal degree, *Adv. Comput. Math*. **4** 389-396 (1995).

16. Z. Wu, Multivariate compactly supported positive definite radial functions, *Adv. Comput. Math*. **4** 283-292 (1995).

17. R. Schaback, Creating surfaces from scattered data using radial basis function. In *Mathematical Methods for Curves and Surfaces*, (Edited by M. Dahlen et al.), pp. 477-496, Vanderbilt University Press, Nashville, (1995).

18. S. M. Wong, Y.C. Hon and M.A. Golberg, Compactly supported radial basis functions for the shallow water equations, *Appl. Math. Comput*. (to appear).

19. C.S. Chen, C.A. Brebbia and H. Power, Boundary element methods using compactly supported radial basis functions, *Commun. Numer. Meth. Engng*. **15** 137-150 (1999).

20. R.L. Hardy, Multiquadratic equations for topography and other irregular surfaces, *J. Geophys. Res*., **176** 1905-1915 (1971).

21. J. Duchon, Interpolation des fonctions de deux variables suivant le principe de la flexion des plaques minces, *RAIRO Analyse Numeriques*, **10** 5-12 (1976).

22. R. Franke, Scattered data interpolation: tests of some methods, *Math. Comput*. **48** 181-200 (1982).

23. Y.C. Hon and X.Z. Mao, An efficient numerical scheme for Burgers' equation, *Appl. Math. Comput*. **95**(1) 37-50 (1998).

24. R.E. Carlson and T.A. Foley, The parameter $R^2$ in multiquadratic interpolation, *Comput. Math. Appl.* **21** 29-42 (1991).

25. P.W. Partridge, C.A. Brebbia and L.W. Wrobel, *The Dual Reciprocity Boundary Element Method*, Comput. Mech. Publ., Southampton, UK (1992).

26. W.H. Press, S.A. Teukolsky, W.T. Vetterling and B.P. Flannery, *Numerical Recipes in Fortran*, Cambridge University Press (1992).

27. C. Canuto, M.Y. Hussaini, A. Quarteroni and T.A. Zang, *Spectral Methods in Fluid Dynamics*, Springer, Berlin (1988).

28. J.M. Ortega, and W.C. Rheinboldt. *Iterative Solution of Nonlinear Equations in Several Variables*, Academic Press, (1970).

29. W. Chen, C. Shu, and W. He, The DQ solution of geometrically nonlinear bending of orthotropic rectangular plates by using Hadamard and SJT product, *Computers & Structures*, 74(1) 65-74 (2000).

30. A.J. Nowak, The multiple reciprocity method of solving transient heat conduction problems. In Boundary Elements XI, Vol. 2, Comput. Mech. Publ., Southamptonn and Springer-Verglag, Berlin and New York (1989).

31. J.T. Katsikadelis and M.S.Nerantzaki, The boundary element method for nonlinear problems, *Engineering Analysis with Boundary Element*, **23** 365-273 (1999).

32. H.C. Yee and R.K. Sweby, Aspects of numerical uncertainties in time marching to steady-state numerical solutions, *AIAA J*. 36 712-723 (1998).

33. S. Zhu, P. Satravaha and X. Lu, Solving linear diffusion equations with the dual reciprocity method in Laplace space. *Engineering Analysis with Boundary Elements*, **13** 1-10 (1994).

34. M. Stehfest, Algorithm 368: numerical inversion of Laplace transform. *Commun. ACM*, **13** 47-49 (1970).

35. G.H. Golub and J.M. Ortega, *Scientific Computing and Differential Equations*, Academic Press, (1992).

36. R. Piltner, Recent development in the Trefftz method for finite element and boundary element application. *Advances in Engineering Software,* **2** 107-115 (1995).

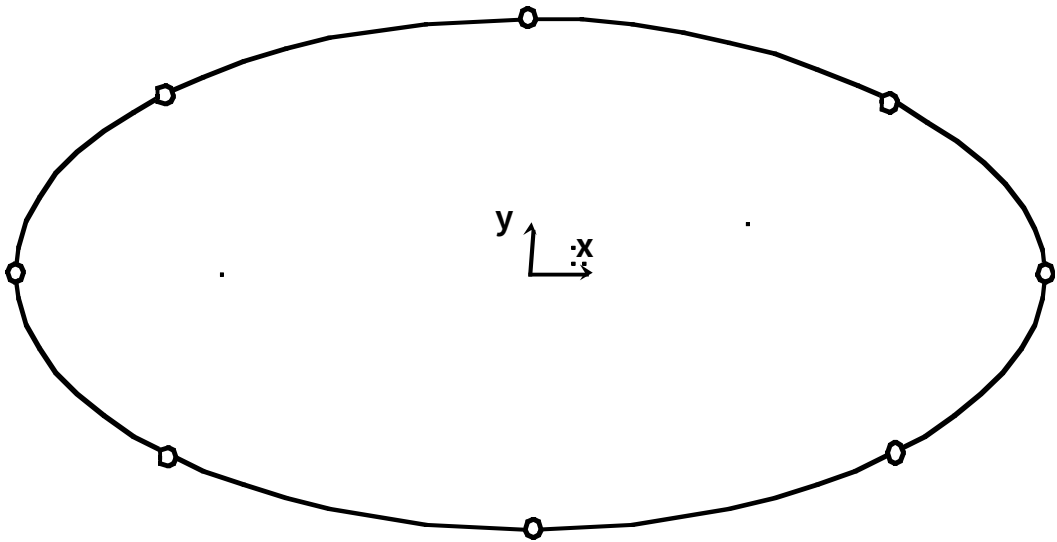

Fig. 1. The BKM boundary knots

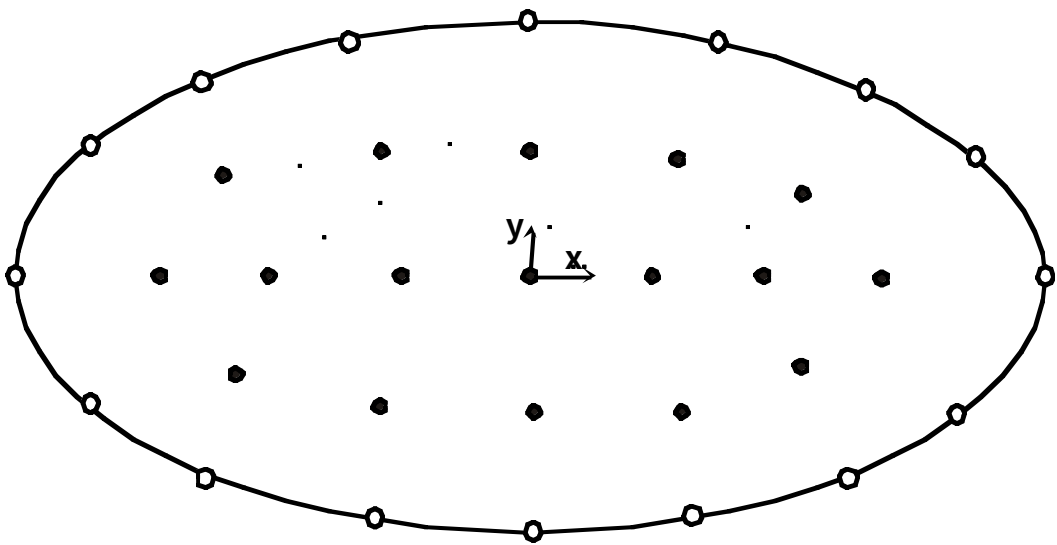

Fig. 2. The DRBEM boundary and internal knots

**Table 1. Results for Helmholtz equation**

| x | y | Exact | DRBEM (33) | BKM (7) | BKM (11) |
|---|---|---|---|---|---|
| 1.5 | 0.0 | 0.997 | 0.994 | 0.999 | 0.997 |
| 1.2 | -0.35 | 0.932 | 0.928 | 0.931 | 0.932 |
| 0.6 | -0.45 | 0.565 | 0.562 | 0.557 | 0.565 |
| 0.0 | 0.0 | 0.0 | 0.0 | 0.0 | 0.0 |
| 0.9 | 0.0 | 0.783 | 0.780 | 0.779 | 0.783 |
| 0.3 | 0.0 | 0.296 | 0.294 | 0.289 | 0.296 |
| 0.0 | 0.0 | 0.0 | 0.0 | 0.0 | 0.0 |

**Table 2. Results for Laplace equation.**

| x | y | Exact | BEM (16) | BKM (3) | BKM (5) |
|---|---|---|---|---|---|
| 1.5 | 0.0 | 1.500 | 1.507 | 1.500 | 1.500 |
| 1.2 | -0.35 | 0.850 | 0.857 | 0.850 | 0.850 |
| 0.6 | -0.45 | 0.150 | 0.154 | 0.150 | 0.150 |
| 0.0 | 0.0 | -0.450 | -0.451 | -0.450 | -0.450 |
| 0.9 | 0.0 | 0.900 | 0.913 | 0.900 | 0.900 |
| 0.3 | 0.0 | 0.300 | 0.304 | 0.300 | 0.300 |
| 0.0 | 0.0 | 0.0 | 0.0 | 0.0 | 0.0 |

**Table 3. Results for $\nabla^2 u = -\partial u/\partial x$**

| x | y | Exact | DRBEM (33) | BKM(15) | BKM (18) |
|---|---|---|---|---|---|
| 1.5 | 0.0 | 0.223 | 0.229 | 0.229 | 0.224 |
| 1.2 | -0.35 | 0.301 | 0.307 | 0301 | 0.305 |
| 0.0 | -0.45 | 1.000 | 1.003 | 1.010 | 1.000 |
| -0.6 | -0.45 | 1.822 | 1.819 | 1.822 | 1.818 |
| -1.5 | 0.0 | 4.482 | 4.489 | 4.484 | 4.477 |
| 0.3 | 0.0 | 0.741 | 0.745 | 0.744 | 0.743 |
| -0.3 | 0.0 | 1.350 | 1.348 | 1.353 | 1.354 |
| 0.0 | 0.0 | 1.000 | 1.002 | 1.003 | 1.004 |

Table 4. Results for $\nabla^2 u = -\partial u/\partial x - \partial u/\partial y$

| x | y | Exact | DRBEM (33) | BKM(15) | BKM (18) |
|---|---|---|---|---|---|
| 1.5 | 0.0 | 1.223 | 1.231 | 1.225 | 1.224 |
| 1.2 | -0.35 | 1.720 | 1.714 | 1.725 | 1.723 |
| 0.0 | -0.45 | 2.568 | 2.557 | 2.546 | 2.551 |
| -0.6 | -0.45 | 3.390 | 3.378 | 3.403 | 3.405 |
| -1.5 | 0.0 | 5.482 | 5.485 | 5.490 | 5.491 |
| 0.3 | 0.0 | 1.741 | 1.731 | 1.729 | 1.731 |
| -0.3 | 0.0 | 2.350 | 2.335 | 2.349 | 2.350 |
| 0.0 | 0.0 | 2.000 | 1.989 | 1.992 | 1.993 |

Table 5. Results for nonlinear equation (44)

| x | y | Exact | BKM (5) | Relative error % |
|---|---|---|---|---|
| 4.5 | 0.0 | 20.25 | 20.34 | -0.44 |
| 4.2 | -0.35 | 17.64 | 17.74 | -0.58 |
| 3.6 | -0.45 | 12.96 | 13.07 | -0.86 |
| 3.0 | -0.45 | 9.00 | 9.09 | -0.98 |
| 2.4 | -0.45 | 5.76 | 5.82 | -0.94 |
| 1.8 | -0.35 | 3.24 | 3.26 | -0.77 |
| 1.5 | 0.0 | 2.25 | 2.25 | -0.18 |
| 3.9 | 0.0 | 15.21 | 15.36 | -0.99 |
| 3.3 | 0.0 | 10.89 | 11.03 | -1.28 |
| 3.0 | 0.0 | 9.00 | 9.12 | -1.32 |
| 2.7 | 0.0 | 7.29 | 7.38 | -1.31 |
| 2.1 | 0.0 | 4.41 | 4.46 | -1.14 |

**Table 6. Results for Burger equation**

| x   | y     | Exact | BKM (5) | Relative error % |
|-----|-------|-------|---------|------------------|
| 4.5 | 0.0   | 0.444 | 0.479   | -7.9             |
| 4.2 | -0.35 | 0.476 | 0.515   | -8.2             |
| 3.6 | -0.45 | 0.555 | 0.585   | -5.4             |
| 3.0 | -0.45 | 0.666 | 0.666   | 0.15             |
| 2.4 | -0.45 | 0.833 | 0.808   | 3.1              |
| 1.8 | -0.35 | 1.111 | 1.089   | 2.0              |
| 1.5 | 0.0   | 1.333 | 1.300   | 2.5              |
| 3.9 | 0.0   | 0.512 | 0.563   | -9.7             |
| 3.3 | 0.0   | 0.606 | 0.632   | -4.2             |
| 3.0 | 0.0   | 0.666 | 0.672   | -7.3             |
| 2.7 | 0.0   | 0.740 | 0.725   | 2.0              |
| 2.1 | 0.0   | 0.952 | 0.918   | 3.6              |

Fig. 1. The BKM boundary-only knots

Fig. 2. The DRBEM boundary and internal knots.